\documentclass[12pt]{amsart}
\usepackage{psfig}
\usepackage{amscd} 
\newtheorem{theo}{Theorem}[section]
\newtheorem{pro}[theo]{Proposition}
\newtheorem{lem}[theo]{Lemma}
\newtheorem{defi}[theo]{Definition}
\newtheorem{remark}[theo]{Remark}
\newtheorem{corol}[theo]{Corollary}

\newtheorem{conjecture}[theo]{Conjecture}
\newtheorem{question}[theo]{Question}

\def \proof {{\bf Proof$\colon$}\ }

\def \Z{{\bf Z}}

\def \N{{\bf N}}
\def \ad {{$K\buildrel n\over \longrightarrow K'$}}

\def \d{$({\mathcal D}, {\mathbf q})\ $}

\begin{document}

\title[Knot adjacency and genus]{Knot adjacency, genus and essential tori}
\author[E. Kalfagianni]{Efstratia Kalfagianni$^1$}
\thanks
{$^1$Supported in part by NSF 
grants DMS-0104000 and DMS-0306995  and by a grant through the Institute for Advanced Study.}
\author[X.-S. Lin]{Xiao-Song Lin$^2$}\thanks{$^2$Supported in part by the
Overseas Youth Cooperation
Research Fund of NSFC  and by NSF
grants DMS-0102231 and DMS-0404511.}

\address[]{Department of Mathematics, Michigan State
University,
East Lansing, MI, 48824}
\email[]{kalfagia@@math.msu.edu}

\address[]{Department of Mathematics, 
University of California,
Riverside, CA, 92521}
\email[]{xl@@math.ucr.edu}

\begin{abstract} A knot $K$ is called $n$-adjacent to another
knot $K'$, if $K$
admits a projection containing $n$ ``generalized crossings"
such that changing any $0<m\leq n$ of them yields a
projection of $K'$.
We apply techniques from the theory of sutured 3-manifolds,
Dehn surgery and the theory of geometric structures of 3-manifolds to
answer the question of
the extent to which non-isotopic knots can be adjacent
to each other. A consequence of
our main result is that if $K$ is $n$-adjacent to $K'$
for all $n\in \N$, then $K$ and $K'$ are isotopic. This
provides a partial verification of the conjecture of V. Vassiliev that
the finite type knot invariants distinguish all knots. We also
show
that if no twist about a crossing circle $L$ of a knot $K$ changes
the isotopy class of $K$, then $L$ bounds a disc
in the complement of $K$. This leads to a characterization of  nugatory
crossings on knots.

\smallskip
\smallskip

\smallskip
\smallskip

\smallskip
\smallskip

\smallskip
\smallskip

\smallskip
\smallskip
{\it AMS classification numbers:} 57M25, 57M27, 57M50.
\smallskip
\smallskip
\smallskip
\smallskip

\end{abstract}

\maketitle

\section{Introduction} A {\it crossing disc} for a knot $K\subset S^3$
is an embedded disc $D\subset S^3$
such
that $K$ intersects ${\rm int}(D)$ twice with
zero algebraic number.  
Let $q\in \Z$.
Performing ${\textstyle {1\over q}}$-surgery on $L_1:=\partial D_1$,
changes $K$ to another knot $K^{'}\subset S^3$.
We will say that 
$K^{'}$ is obtained from $K$
by a generalized crossing change of order $q$ (see Figure 1).

An $n$-{\it collection} for a knot $K$ is a pair
$({\mathcal D}, {\mathbf q})$ , such that:
\vskip .04in

i) ${\mathcal D}:=\{  D_1,\ldots, D_{n}\}$
is a set of
disjoint
crossing discs for
$K$;
\vskip .04in

ii)  ${\mathbf q}:=\{ {\textstyle {1\over q_1}},\ldots,
{\textstyle {1\over q_n}}\}$,
with $q_i \in \Z-\{0\}$ ;
\vskip .04in

iii)
the knots
$ L_1,\ldots, L_n$
are labeled by
${\textstyle {1\over q_1}},\ldots,
{\textstyle {1\over q_n}}$,
respectively. Here, $L_i:=\partial D_i$.
The link $L:=\cup_{i=1}^n L_i$ is called the {\it crossing link} 
associated to \d.
\vskip .08in

Given a knot $K$ and an $n$-collection \d,
for $j=1, \ldots, n$, let $i_j \in \{0, \ 1\}$ and
$${\mathbf i}:=(i_1,\ldots, i_{n}).$$
We denote ${\mathbf 0}=(0,\ldots,0)$ and ${\mathbf 1}=(1,\ldots,1)$.

For every ${\mathbf i}$,
we will denote by $K({\mathbf i})$
the knot obtained from $K$ by a surgery modification
of order $q_i$ (resp. $0$),
along each $L_j$ for
which $i_j=1$ (resp. $i_j=0$). 

\begin{defi} \label{defi:surgery} We will say that $K$
is $n$-adjacent to $K'$ 
if there exists 
an { $n$-collection} \d  for $K$,
such that
the knot  $K({\mathbf i})$
is isotopic to $K'$
for every ${\mathbf i}\neq {\mathbf 0}$. We will write \ad\  and we will say that
\d transforms $K$ to $K'$. 
\end{defi}

\begin{figure}[ht]
\centerline{\psfig{figure=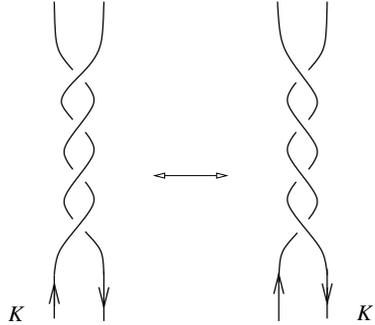,height=1.7in, clip=}}
\caption{The knots $K$ and $K'$ differ by a generalized crossing change
of order $q=-4$.}
\end{figure}

Our main result is the following:

\begin{theo}\label{theo:main} Suppose that $K$ and $K'$
are non-isotopic knots. 
There exists a
constant $C(K,\,K')$ such that
if \ad, 
then 
$n\leq C(K,\,K')$.
\end{theo}

The quantity $C(K,\,K')$ can be expressed in terms  of computable invariants
of the knots $K$ and $K'$.
Let
$g(K)$ and $g(K')$ denote the genera of $K$ and $K'$, respectively
and let
$g:={\rm max}\,\{\,g(K),\,g(K')\,\}$.
The  constant
$C(K,\,K')$ encodes
information about the relative size of
$g(K)$, $g(K')$
and the behavior of the satellite structures of $K$ and
$K'$ under the Dehn surgeries
imposed by knot adjacency.
In many cases $C(K,\,K')$
can be made explicit. For example,
when $g(K)>g(K')$
we have
$C(K,\,K')=6g-3$. Thus, in this case,
Theorem \ref{theo:main} can be restated as follows:

\begin{theo} \label{theo:lower}
Suppose that $K,K'$ are knots with
$g(K)>g(K')$.
If  \ad,  
then $n\leq 6g(K)-3$.
\end{theo}

In the case that $K'$ is the trivial knot
Theorem \ref{theo:lower} was proven by
H. Howards
and J. Luecke in \cite{kn:hl}.

A crossing of a knot $K$, with crossing disc $D$, is called
{\em nugatory} iff $\partial D$ bounds a disc that is disjoint from $K$.
The techniques used in the proof of Theorem \ref{theo:main}
have applications to
the question of whether a crossing change that doesn't change the isotopy
class of the underlying knot
is nugatory (Problem 1.58, \cite{kn:kirby}). As a corollary of the proof of
Theorem \ref{theo:main}
we obtain the following characterization of nugatory crossings:

\begin{corol} \label{corol:nstrong}
For a crossing disc $D$ of a knot $K$ let $K(r)$ denote the knot obtained
by a twist of order $r$  along $D$. The crossing is nugatory
if and only if
$K(r)$ is isotopic to $K$ for all $r\in \Z$. 
\end{corol}

Definition \ref{defi:surgery} is equivalent to the definition of
$n$-adjacency given in the abstract of this paper or in \cite{kn:kl}. With this reformulation,
it follows that $n$-adjacency implies $n$-similarity
in the sense of \cite{kn:Oh}, which in turn, as shown in \cite{kn:NS}, implies
$n$-equivalence. Gussarov showed that
 two knots are $n$-equivalent
precisely when all of their finite type invariants of orders
$<n$ are the same. 
Vassiliev (\cite{kn:vas}) has conjectured that if two oriented knots
have all of their finite type invariants the same then
they are isotopic.
In the light Gussarov's result, Vassiliev's
conjecture
can be reformulated as follows:

\begin{conjecture}
\label{conjecture:vas}Suppose that $K$ and $K'$ are knots that are
$n$-equivalent for all $n\in \N$. Then $K$ is isotopic to $K'$.
\end{conjecture}

To that respect, Theorem \ref{theo:main} implies the following corollary
that provides a partial verification to Vassiliev's conjecture:

\begin{corol}\label{corol:afirm}
If \ad, for all $n\in \N$, then $K$ and $K'$ are isotopic.
\end{corol}
We now describe the contents of the paper and the idea of the proof
of the main theorem.
Let $K$ be a knot and let \d be a $n$-collection
with associated crossing link $L$. 
Since the linking number
of
$K$ and every component of $L$ is zero,
$K$ bounds a Seifert surface in the complement
of $L$. Thus, we can define the genus of $K$ in the complement
of $L$, say $g_L^n(K)$.
In Section two we 
study the question of the extent to which a Seifert surface of
$K$ that is of minimal genus in the complement
of $L$ remains of minimal genus under various 
surgery modifications along the components of $L$.
Using a result of of Gabai (\cite{kn:ga})
we show  that if \ad, and \d is an $n$-collection
that transfers $K$ to $K'$ then
$g_L^n(K)=g:={\rm max}\,\{\,g(K),\, g(K')\,\}$,
where $g(K), g(K')$ denotes the genus of $K$, $K'$
respectively. This is done in Theorem \ref{theo:genus}.

In Section three, we prove
Theorem \ref{theo:lower}.
In Section four, we finish the proof of Theorem
\ref{theo:main}:
We begin by defining a notion
of $m$-adjacency between knots $K, K'$  with respect to an one component
crossing link $L_1$ of $K$ (see Definition \ref{defi:curve}).
To describe our approach in more detail,
set $N:=S^3\setminus \eta(K\cup L_1)$,
and
let $\tau(N)$  denote the number of
disjoint, pairwise non-parallel, essential
embedded tori in $N$.
We employ results of Cooper and Lackenby (\cite{kn:cla}), Gordon
(\cite{kn:go}) and  McCullough (\cite{kn:appen})
and an 
induction argument on $ \tau(N)$
to show the following: Given knots $K, K'$,
there exists a constant $b(K,\ K')\in {\bf N}$
such that if $K$ is $m$-adjacent to $K'$
with respect to a crossing link $L_1$
then either $m\leq b(K,\ K')$ or
$L_1$ bounds an embedded disc in the
complement of $K$. This is done in Theorem \ref{theo:last}.
Theorem
\ref{theo:genus} implies that if
\ad\,  
and $n>m(6g-3)$,
then an $n$-collection
that transforms $K$ to $K'$ gives rise to a
crossing link $L_1$
such that $K$ is $m$-adjacent to $K'$
with respect to $L_1$.
Combining this with Theorem \ref{theo:last} yields Theorem \ref{theo:main}.

In Section five, we present some applications
of the results of Section four and the methods used in their proofs.
Also,
for every $n\in \N$, we construct examples
of non-isotopic knots $K,K'$ such that \ad.

Throughout the entire paper we work in the PL
or the smooth category.
In \cite{kn:k1},
the techniques of this paper are refined and used to study adjacency to fibered knots and the problem
of nugatory crossings in fibered knots.
In \cite{kn:kl1} the results of this paper
are used to obtain criteria for detecting non-fibered knots and
for detecting the non-existence of symplectic
structures on certain 4-manifolds. Further applications
include constructions of 3-manifolds that are indistinguishable by certain Cochran-Melvin finite type invariants (\cite{kn:kl1}),
and constructions of hyperbolic knots with trivial Alexander polynomial
and arbitrarily large volume (\cite{kn:k}).
\vskip 0.04in

{\bf Acknowledgment.} We thank Tao Li,
Katura Miyazaki and Ying-Qing Wu for their interest in this work and
for their helpful comments on  an earlier version of the paper.
We thank Darryl McCullough for his comments and for proving
a result about homeomorphisms of 3-manifolds (\cite{kn:appen}) that is needed for the proof the
main result of this paper.
We also thank Ian Agol, Steve Bleiler, Dave Gabai, Marc Lackenby,
Marty Scharlemann and Oleg Viro for useful conversations or correspondence.
Finally, we thank the referee for a very thoughtful and careful review that has lead 
to a significant improvement of the exposition in this paper.

\smallskip
\smallskip

\section{Taut surfaces, knot genus and multiple crossing changes}
Let $K$ be a knot and \d an $n$-collection for $K$
with associated crossing link
$L$.
Since the linking number
of
$K$ and every component of $L$ is zero,
$K$ bounds a Seifert surface $S$ in the complement
of $L$. Define
$$g_n^L(K):={\rm min}\,\{\, {\rm genus (S)}\,|\,S\, {\rm a\ Seifert \ surface\  of \,}K\, {\rm as\ above}\,\}.$$
Our main result in this section is the following:

\begin{theo} \label{theo:genus} Suppose that \ad, for some $n\geq 1$. Let \d
be an $n$-collection that transforms $K$ to $K'$
with
associated crossing link $L$. 
We have
$$g_n^L(K)={\rm max}\,\{\, g(K),g(K')\,\}.$$
In particular, $g_n^L(K)$ is independent of $L$ and $n$.
\end{theo}

Before we proceed with the proof of Theorem \ref{theo:genus} we need some preparation:
For a link ${\bar L}\subset S^3$ we will use 
$\eta({\bar L})$ to denote a regular neighborhood of $\bar L$.
For a knot $K\subset {S}^3$ and an $n$-collection \d,
let $$M_{L}:={S}^3 \setminus \eta(K\cup L),$$
where $L$ is the crossing link associated to \d.

\begin{lem} \label{lem:irreducible}Suppose that $K$, $K'$ are
knots such that \ad, for some $n\geq 1$. Let \d be an $n$-collection
that transforms $K$ to $K'$.
If $M_L$ is reducible then a component of $L$ bounds an embedded disc
in the complement of $K$. Thus, in particular,
$K$ is isotopic to $K'$.
\end{lem}

\proof Let $\Sigma$ be an essential 2-sphere in $M_L$.
Assume that $\Sigma$ has been isotoped so that the intersection
$I:=\Sigma \cap ({\cup_{i=1}^nD_i})$ is minimal.
Notice that we must have $I\neq \emptyset$
since otherwise $\Sigma$ would bound
a 3-ball in $M_L$. Let $c\in (\Sigma\cap D_i)$ denote a component of $I$ that is innermost
on $\Sigma$; that is $c$ bounds a disc $E\subset \Sigma$
such that ${\rm int}(E)\cap ({\cup_{i=1}^n D_i})=\emptyset$.
Since $\Sigma$ is separating in $M_L$, $E$ can't contain just one point of $K\cap D_i$. $E$ can't be disjoint from $K$ or $c$ could be removed by isotopy. Hence $E$ contains both points of $K\cap D_i$ and so $c=\partial E$ is parallel to $\partial D_i$ in $D_i\setminus K$.
It follows that $L_i$ bounds an embedded disc
in the complement of $K$. Since
${\textstyle {1\over q_i}}$-surgery on $L_i$
turns $K$ into $K'$, we conclude that 
$K$ is isotopic to $K'$. \qed

\smallskip

\vskip .03in

To continue we recall the following definition:
\vskip .03in

\begin{defi}\label{defi:thurston}{\rm (\cite{kn:th})} Let $M$ be
a compact, oriented 3-manifold with boundary $\partial M$.
For a compact, connected, oriented surface
$(S, \ \partial S) \subset (M, \ \partial M)$,
the complexity $\chi^{-}(S)$ is defined by 

$\chi^{-}(S):={\rm max}\, \{\,0,\, -\chi(S)\,\}$,
where $\chi(S)$ denotes the Euler
characteristic of $S$. If $S$ is disconnected
then $\chi^{-}(S)$ is defined to be the sum of the complexities of all the components of $S$.
Let $\eta(\partial S)$ denote a regular neighborhood
of $\partial S$ in
$\partial M$.
The Thurston norm $x(z)$ of a homology class
$z\in H_2(M, \ \eta(\partial S))$
is the
minimal
complexity over all oriented, embedded surfaces
representing $z$. 
The surface $S$ is 
called {\it taut} if
it is incompressible
and 
we have $x([S, \ \partial S])=\chi^{-}(S)$; that is
$S$ is norm-minimizing.
\end{defi}
\vskip 0.04in

We will need the following lemma the proof of which follows  from
the definitions:

\begin{lem} \label{lem:taut} Let \d be an $n$-collection
for a knot $K$ with
associated crossing link $L$ and $M_L:=S^3\setminus\eta(K\cup L)$.
A compact, connected, oriented surface 
$(S, \ \partial S )\subset (M_L, \ \partial \eta(K))$,
such that $\partial S=K$, is taut
if and only if among all Seifert surfaces
of $K$ in the complement of $L$,
$S$ has the minimal genus.
\end{lem}
\smallskip

To continue, we need to introduce some more notation.
For ${\mathbf i}$ as before the statement of Definition
\ref{defi:surgery}, let
$M_L({\mathbf i})$ denote
the 3-manifold obtained from $M_L$
by performing Dehn filling on $\partial M_L$ as follows:
The slope of the filling for the components
$\partial \eta(L_j)$ for
which $i_j=1$ (resp. $i_j=0$) is 
${\textstyle {1\over q_j}}$ (resp. $\infty:= {\textstyle {1\over 0}}$).
Clearly we have $M_L({\mathbf i})=S^3\setminus \eta(K({\mathbf i}))$,
where $K({\mathbf i})$ is as in Definition \ref{defi:surgery}.
Also let $M^{+}_L({\mathbf i})$
(resp. $M^{-}_L({\mathbf i})$) denote
the 3-manifold obtained from $M_L$
by only performing Dehn filling
with slope
${\textstyle {1\over q_j}}$
(resp. $\infty$)
on 
the components
$\partial \eta(L_j)$ for
which $i_j=1$. 

\begin{lem} \label{lem:key}
Let \d be an $n$-collection for a knot $K$
such that
$M_L$ is irreducible. 
Let $(S, \ \partial S )\subset (M_L, \ \partial \eta(K))$
be an oriented surface with
$\partial S=K$
that is taut. For $j=1, \ldots, n$,
define ${\mathbf i}_j:=(0, \ldots, 0,  1, 0, \ldots , 0)$
where the unique entry 1 appears at the $j$-th place.
Then,
at least one of $M^{+}_L({\mathbf i}_j)$, $M^{-}_L({\mathbf i}_j)$
is irreducible and $S$ remains taut in that 3-manifold.
\end{lem}
\proof The proof uses a result of \cite{kn:ga}
in the spirit of \cite{kn:st}:
For $j\in \{1,\ldots, n\}$ set  $M^{+}:=M^{+}_L({\mathbf i}_j)$ and
$M^-:= M^{-}_L({\mathbf i}_j)$. Also set
$L^j:= L \setminus L_j$ and $T_j:=\partial\eta(L_j)$.
We distinguish two cases:
\smallskip

 {\it Case 1:} Suppose that every embedded torus that is incompressible
in $M_L$
and it separates $L^j\cup S$ from $L_j$,
is parallel to $T_j$.
Then, $M_L$ is $S_{L_j}$-atoroidal (see
Definition 1.6 of
\cite{kn:ga}). By Corollary 2.4 of \cite{kn:ga},
there is at most one Dehn filling
along
$T_j$ that yields a
3-manifold  which is either reducible or in which
$S$ doesn't remain taut.
Thus the desired conclusion follows. 
\smallskip

 {\it Case 2:} There exists  an
embedded
torus
$T\subset M_L$ such that
i) $T$ is incompressible in $M_L$;
ii) $T$ separates $L^j\cup S$ from $L_j$; and
iii) $T$
is not parallel to $T_j$.
In $S^3$, $T$ bounds a
solid torus $V$, with $\partial V=T$. 
Suppose, for a moment, that $L_j$ lies in ${\rm int}(V)$
and $L^j\cup S$ lies in $S^3\setminus V$. 
If $V$ is knotted in $S^3$ then, since
$L_j$ is unknotted, $L_j$ is homotopically inessential in $V$. But then
$T$  compresses
in $V$ and thus in $M_L$; a contradiction. If $V$ is unknotted in $S^3$ then the longitude of $V$ bounds a
disc $E$
in $S^3\setminus V$. Since $S$ is disjoint from
$T$, $K$ intersects $E$ at least twice.
On the other hand, since $T$ is incompressible in $M_L$
and $K$ intersects $D_j$ twice,
$L_j$ is isotopic to the core of $V$. 
Hence, $T$ is parallel to $T_j$ in $M_L$; a contradiction.
Hence $L^j\cup S$ lies in ${\rm int}(V)$ while 
$L_j$ lies in $S^3\setminus V$. We will show that
$M^{+}$, $M^{-}$ are irreducible and that
$S$ remains taut in both of these 3-manifolds.

Among all tori in $M_L$ that have properties (i)-(iii)
stated above, choose $T$ to be one that minimizes
$|T\cap D_j|$. Then,  that $D_j\cap T$
consists of a single curve which bounds a
disc $D^*\subset {\rm int}( D_j)$, such that
$(K\cap D_j)\subset {\rm int}( D^*)$
and $D^*$ is a meridian disc of $V$. See Figure 2 below.
\begin{figure}[ht]
\centerline{\psfig{figure=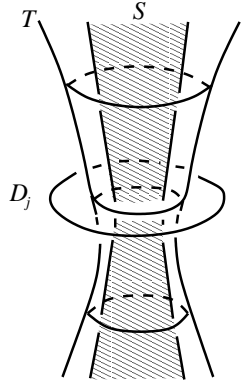,height=2in, clip=}}
\caption{The intersection of $T$ and $S$ with $D_j$.}
\end{figure}
Since $T$ is not parallel to $T_j$,
$V$ must be knotted. 
For $r\in {\bf Z}$, 
let $M(r)$ denote the 3-manifold obtained from $M_L$
by performing Dehn filling along $\partial \eta(L_j)$
with slope 
$\textstyle{1\over {r}}$. 
Since the core of $V$ intersects $D_j$
once, the Dehn filling doesn't unknot $V$ and
$T=\partial V$ remains incompressible in $M(r)\setminus V$.
On the other hand, $T$ is incompressible in $V\setminus(K\cup L^j)$
by definition. Notice that both $M(r)\setminus V$ and $V\setminus(K\cup L^j)$ are irreducible and 
$$M(r)=(M(r)\setminus V)\,\bigcup_T\,(V\setminus(K\cup L^j)).$$ 
We conclude that 
$T$ remains incompressible in $M(r)$ and $M(r)$ is irreducible.
In particular
$M^{+}$ and $M^{-}$ are both irreducible.

Next we show that $S$ remains taut in $M^{+}$ and $M^{-}$.
By
Lemma \ref{lem:taut}, we must show that
$S$ is a minimal genus surface for 
$K$ in $M^{+}$ and in $M^-$.
To that end, let $S_1$
be a minimal genus surface for $K$ in $M^{+}$ or in $M^-$.
We may isotope so that $S_1\cap T$ is
a collection of parallel essential curves on $T$.
Since the linking number of $K$ and $L_j$ is zero,
$S_1\cap T$ is homologically trivial in
$T$. Thus, we may attach annuli along the components of 
$S_1\cap T$ and then isotope off $T$ in ${\rm int}(V)$,
to obtain
a Seifert surface $S_1'$ for $K$ that is disjoint from $L_j$.
Thus $S_1'$  is a surface in the complement of $L$.
Since $T$ is incompressible, no component of $S_1\setminus V$ is a disc.
Thus,
${\rm genus (S_1')}\leq {\rm genus (S_1)}$.
On the other hand, by definition of $S$,
${\rm genus (S)}\leq {\rm genus (S_1')}$
and thus ${\rm genus (S)}\leq {\rm genus (S_1)}$.  \qed
\smallskip

\begin{lem} \label{lem:key2}
Let \d be an $n$-collection for a knot $K$
such that
$M_L$ is irreducible. 
Let $(S, \ \partial S )\subset (M_L, \ \partial \eta(K))$
be an oriented surface with
$\partial S=K$
that is taut. There exists at least one 
sequence 
${\mathbf i}:=(i_{1},\ldots i_{n})$, with $i_{j}\in \{1, 0\}$,
such that $S$ remains taut in 
$M_{L}({\mathbf i})$. Thus we have, $g(K({\mathbf i}))={\rm genus}(S)$ .
\end{lem}

\proof The proof is by induction on $n$.
For $n=1$, the conclusion follows from Lemma \ref{lem:key}.
Suppose inductively that for every $m<n$
and every $m$-collection
$({\mathcal D}_1, {\mathbf q}_1)$ of a knot $K_1$ such that
$M_{{L}^1}$ is irreducible, the conclusion of the lemma
is true. Here,
$L^1$ denotes
the crossing link associated to  
${\mathcal D}_1$ and 
$M_{L^1}:= S^3 \setminus \eta(K_1\cup L^1)$.

Now let $K$, \d  and $S$ be as in the statement of the lemma and
let ${\mathbf i}_1:=(1, 0,\ldots, 0)$. By Lemma \ref{lem:key} at least one of
$M_L^{\pm}({\mathbf i}_1)$,
say $M_L^{-}({\mathbf i}_1)$, is irreducible and $S$ remains
taut in that 3-manifold.
Let 
$${\mathcal D}_1:=\{D_2, \ldots, D_n\} \ {\rm and } \ 
{\mathbf q}_1:=\{q_2, \ldots, q_n\}.$$
Let ${{L}^1}:=L\setminus L_1$ and let $K_1$ denote the image
of $K$ in $M_L^{-}({\mathbf i}_1)$.
Clearly, $M_{L^1}=M_L^{-}({\mathbf i}_1)$ and thus $M_{L^1}$ is irreducible.
By the induction hypothesis, applied to $K_1$ and the 
$(n-1)$-collection $({{\mathcal D}_1}, {\mathbf q}_1)$,
it follows that
there is at least one 
sequence 
${\mathbf i}_0:=(i_{02},\ldots i_{0n})$, with $i_{0j}\in \{1, 0\}$,
such that $S$ remains taut in $M_{L^1}({\mathbf i}_0)$.
Since $M_{L^1}({\mathbf i}_0)=M_L({\mathbf i})$,
where ${\mathbf i}:=(0, i_{02},\ldots i_{0n})$, the desired conclusion
follows. \qed
\medskip
\smallskip

\proof [Proof of Theorem \ref{theo:genus}]
Let \ad, $L$ and $M_L$ be as in the statement
of the theorem.
Let $S$ be a Seifert surface
for $K$ in the complement of $L$
such that ${\rm genus}(S)=g_n^L(K)$. 
First, assume that $M_L$ is irreducible.
By Lemma \ref{lem:taut}, $S$ gives rise to a surface
$(S, \ \partial S)\subset (M_L, \ \eta(\partial S))$
that is {\it taut}. By Lemma \ref{lem:key2},
there exists at least one
sequence 
${\mathbf i}:=(i_{1},\ldots i_{n})$, with $i_{j}\in \{1, 0\}$,
such that $S$ remains taut in 
$M_{L}({\mathbf i})$. There are three cases to consider:

(1) $g(K)>g(K')$,

(2) $g(K)<g(K')$,

(3) $g(K)=g(K')$.

In case (1), for every ${\mathbf i}\neq{\mathbf 0}$, we have
$$g(K')=g(K({\mathbf i}))<g(K)\leq{\rm genus}(S).$$
Therefore $S$ doesn't remain taut in
$M_{L}({\mathbf i})=S^3\setminus \eta(K({\mathbf i}))$.
Hence
$S$ must
remain taut in $M_L({\mathbf 0})=S^3\setminus \eta(K)$
and we have $g_n^L(K)=g(K)$.
In case (2), notice that we have a $n$-collection $({\mathcal D}',{\mathbf q}')$
for $K'$ where ${\mathcal D}'={\mathcal D}$ and ${\mathbf q}'=-\mathbf q$,
such that $K'({\mathbf i})=K'$ for all ${\mathbf i}\neq{\mathbf 1}$
and $K'({\mathbf 1})=K$. So we may argue similarly as in case (1) that
$g_n^L(K)=g(K')$. In fact, in case (2), $S$ must remain taut in
$M_L({\mathbf i})$ for all ${\mathbf i}\neq {\mathbf 0}$.
Finally in case (3),
$S$ remains taut in
$M_L({\mathbf i})$ for all ${\mathbf i}$, and it follows that
$g_n^L(K)=g(K')=g(K)$.

Suppose, now, that $M_L$
is reducible. 
By Lemma
\ref{lem:irreducible}, there is 
at least one component of $L$ that bounds
an embedded disc in the complement of $K$. Let $L^1$ denote the union
of the components of $L$ that bound disjoint
discs in the complement of $K$
and let $L^2:=L\setminus L^1$.
We may isotope $S$ so that it is disjoint
from the discs bounded by the components of
$L^1$. Now $S$ can be viewed as taut surface
in $M_{L^2}:=S^3\setminus \eta(K\cup L^1)$.
If $L^2=\emptyset$, the
conclusion is clearly true. Otherwise $M_{L^2}$ is irreducible
and the argument described above applies. \qed
\smallskip

\section{Genus reducing n-collections}
The purpose of this section is to prove Theorem \ref{theo:main}
in the case that
$g(K)>g(K')$.
The argument 
is essentially that in
the proof of  the main result of
\cite{kn:hl}.
\smallskip

\proof [Proof of Theorem \ref{theo:lower}]
Let $K$, $K'$ be 
as in the statement of the theorem.
Let \d be an $n$-collection that transforms $K$ to $K'$
with associated crossing link  $L$.  Let $S$ be a Seifert surface for $K$ that is of minimum
genus among all surfaces bounded by $K$ in 
the complement of $L$. By Theorem
\ref{theo:genus} we have ${\rm genus}(S)=g(K)$.
Since $S$ is incompressible, after an isotopy, we can arrange
so that for $i=1, \ldots, n$, each closed component of
$S\cap {\rm int}(D_i)$
is essential in $D_i\setminus K$ and thus
parallel to $L_i=\partial D_i$ on $D_i$. 
Then , after an isotopy of $L_i$ in the
complement of $K$,
we may assume that
 $S\cap {\rm int}(D_i)$
consists of a single properly embedded arc
$(\alpha_i, \partial \alpha_i) \subset (S,\partial S)$ (see Figure 3).
Notice that $\alpha_i$ is essential on $S$.
For, otherwise, 
$D_i$ would bound a disc in the complement of $K$
and thus the genus of $K$ could
not be lowered by surgery on $L_i$.

{\vspace{.03in}}
\begin{figure}[ht]
\centerline{\psfig{figure=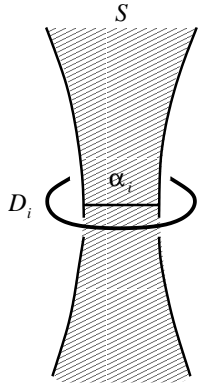,height=2in, clip=}}
\caption{The intersection of $S$ with ${\rm int}(D_i)$.}
\end{figure}

We claim that no two of the arcs $\alpha_1,\ldots \alpha_n,$
can be parallel on $S$.
For, suppose 
on the contrary, that the arcs $\alpha_i:= {\rm int}(D_i) \cap S$
and $\alpha_j:= {\rm int}(D_j) \cap S$ are parallel on $S$.
Then the crossing circles 
$L_{i}$ and $L_{j}$
cobound an embedded annulus that is disjoint from $K$.
Let $$M:=S^3\setminus \eta(K\cup L_i)\ {\rm and} \ 
M_1:=S^3\setminus \eta(K\cup L_i\cup L_j).$$
For $r,s\in {\bf Z}$ let
$M(r)$ (resp. $M_1(r, s)$)
denote the 3-manifold obtained from
$M$ (resp. $M_1$) by filling in 
$ \partial \eta(L_{i})$ (resp. $\partial \eta(L_{i}\cup L_{j})$) 
with slope ${\textstyle {1\over {r}}}$
(resp. slopes ${\textstyle {1\over {r}}}, {\textstyle {1\over {s}}}$).
By assumption,
$S$ doesn't remain taut in any of 
$M(q_i)$, $M_1(q_i,q_j)$.
Since $L_{i}, L_{j}$ are coannular
we see that $M_1(q_i,q_j)=M(q_i+q_j)$.
Notice that $q_i+q_j\neq q_i$ since otherwise 
we would conclude that a twist of order $q_j$
along $L_j$ cannot reduce the genus of $K$. Hence we would have
two distinct Dehn fillings of $M$ along $\partial \eta(L_{i})$ under which $S$
doesn't remain taut, contradicting Corollary 2.4 of \cite{kn:ga}.
Therefore, we conclude that no two of the arcs $\alpha_1,\ldots \alpha_n,$
can be parallel on $S$.
Now the conclusion follows
since a Seifert surface of genus $g$ 
contains
$6g-3$ essential arcs no pair of which is parallel. \qed
\smallskip

\section{Knot adjacency and essential tori} In this section
we will complete the proof
of Theorem \ref{theo:main}. For this we need
to study the case of $n$-adjacent knots  \ad\
in the special situation where all the crossing changes from $K$ to $K'$ are supported
on a single crossing circle of $K$. Using Theorem \ref{theo:genus},
we will see that the general case is reduced to this special one.
\smallskip

\subsection{Knot adjacency with respect to a crossing circle}
We begin with the following definition that provides a refined version of knot adjacency:

\begin{defi} \label{defi:curve} Let $K$, $K'$ be knots and let
$D_1$ be a crossing disc for $K$. We
will say that $K$ is $m$-adjacent to $K'$
with respect to the crossing circle $L_1:=\partial D_1$, if there exist non-zero
integers $s_1, \ldots, s_m$ such that the following is true:
For every $\emptyset \neq J\subset \{1,\ldots, m\}$,
the knot obtained from $K$ by a surgery modification of order
$s_J:= \sum_{j\in J} s_j$ along $L_1$
is isotopic to $K'$.
We will write $K\buildrel {(m, L_1)} \over \longrightarrow K'$.
\end{defi}

Suppose that $K\buildrel {(m,  L_1)} \over \longrightarrow K'$
and consider the $m$-collection obtained 
by taking $m$ parallel copies
of $D_1$ and labeling the $i$-th copy
of $L_1$ by ${\textstyle {1\over s_i}}$.
As it follows immediately from the
definitions, this $m$-collection
transforms $K$ to
$K'$
in the sense of Definition \ref{defi:surgery};
thus $K\buildrel {m} \over \longrightarrow K'$.
The following lemma provides a converse
statement that
is needed for the proof of Theorem \ref{theo:main}:

\begin{lem} \label{lem:relation} Let $K, K'$ be knots and set
$g:={\rm max}\,\{\,g(K), g(K')\,\}$. Suppose that \ad.
If $n>m(6g-3)$ for some $m>0$, then
there exists a crossing link $L_1$ for $K$ 
such that $K\buildrel {(m+1,  L_1)} \over \longrightarrow K'$.
\end{lem}
\proof
Let \d be an $n$-collection that transforms $K$ to $K'$
and let $L$ denote the associated crossing link. 
Let $S$ be a Seifert surface for $K$ that is of minimal
genus among all surfaces bounded by $K$ in 
the complement of $L$.  
Isotope so that, for $i=1, \ldots, n$, the intersection $S\cap {\rm int}(D_i)$
is an arc $\alpha_i$ that is properly embedded and essential on $S$.
By Theorem \ref{theo:genus}, we have
${\rm genus}(S)=g$. Since $n> m(6g-3) $,
the set $\{\alpha_i\ |\ i=1,\ldots, n\}$
contains at least $m+1$ arcs
that are parallel on $S$. Suppose, without loss of generality,
that these are the arcs $\alpha_i$, $i=1,\ldots, m+1$.
It follows that the components $L_1,\ldots, L_{m+1}$
of $L$ are isotopic in the complement of $K$;
thus any surgery along any of these components can be realized as surgery on $L_1$.
It now follows from Definitions \ref{defi:surgery}
and \ref{defi:curve} that $K\buildrel {(m+1,  L_1)} \over \longrightarrow K'$. \qed
\smallskip

The main ingredient needed to complete the proof of Theorem
\ref{theo:main} is provided by the following theorem:

\begin{theo} \label{theo:last} Given knots $K$, $K'$,
there exists a constant $b(K, \ K')\in {\bf N}$, that depends only on $K$ and $K'$,
such that 
if $L_1$ is a
crossing circle of $K$
and
$K\buildrel {(m, L_1)} \over \longrightarrow K'$, then either
$m\leq b(K, \ K')$ or
$L_1$ bounds an embedded disc in the complement
of $K$.
\end{theo}
\smallskip

\proof {[Proof of Theorem \ref{theo:main} assuming Theorem \ref{theo:last}]}
Suppose that $K, K'$ are non-isotopic knots with
\ad. 
If $g(K)>g(K')$ the conclusion follows from Theorem \ref{theo:lower} by
simply taking $C(K,\,K'):=6g-3$.
In general, let $C(K,\  K'):=b(K,\ K')\  (6g-3)$,
where $b:=b(K,\ K')$ is 
the constant of Theorem \ref{theo:last}.
We claim that 
we must have $n\leq C(K,\  K')$.
Suppose, on the contrary, that
$n>C(K,\  K')$.
By Lemma \ref{lem:relation},
there exists a crossing circle $L_1$
for $K$, such that
$K\buildrel {(b+1, L_1)} \over \longrightarrow K'$.
By Theorem \ref{theo:last}, 
$L_1$ bounds an embedded disc in the complement
of $K$. But this implies
that $K$ is isotopic to $K'$ contrary to our assumption. \qed
\smallskip

\smallskip

The rest of this section will be devoted to the proof of Theorem \ref{theo:last}.
For that we need to study whether the complement of $K\cup L_1$ contains essential tori and how these tori
behave under the crossing changes from $K$ to $K'$. Given $K, K'$ and $L_1$ such that
$K\buildrel {(m, L_1)} \over \longrightarrow K'$, set
$N:=S^3\setminus \eta (K\cup L_1)$ and $N':=S^3\setminus \eta(K')$.
By assumption, $N'$ is obtained by Dehn filling 
along the torus
$T_1:=\partial \eta(L_1)$.
If $N$ is reducible, Lemma \ref{lem:irreducible} implies that
$L_1$ bounds a disc in the complement of $K$;
thus Theorem \ref{theo:last} holds.
For irreducible $N$, as it turns out,
there are three basic cases to consider:
\vskip 0.04in

(a) $K'$ is a composite knot.

(b) $N$ is atoroidal.

(c) $N$ is toroidal and $K'$ is not a composite knot.
\vskip 0.04in

By Thurston (\cite{kn:th1}), if $N$ is atoroidal
then it is either hyperbolic (it admits a complete hyperbolic metric of finite volume)
or it is a  Seifert fibered space. To handle the hyperbolic case we use a result
of Cooper and Lackenby (\cite{kn:cla}). The Seifert fibered spaces that occur
are known to be very special and  this case is handled by a case-by-case analysis.
Case (c)  is handled by induction on the number of essential tori
contained in $N$. To set up this induction  one needs to study the behavior of these essential tori under the Dehn fillings from $N$ to $N'$. In particular, one needs to know the circumstances under which
these Dehn fillings {\em create} essential tori in $N'$. For this step, we employ a result of Gordon
(\cite{kn:go}).
\smallskip

\subsection{Composite knots} Here we examine the
circumstances under which a knot $K$ is
$n$-adjacent to a composite knot $K'$. 
We will need the following
theorem.

\begin{theo} \label{theo:composite} \ {\rm  (Torisu,\ \cite{kn:to})} \ Let $K':=K'_1\# K'_2$ be a composite knot and $K''$
a knot that is obtained from $K'$ by a generalized crossing change with corresponding
crossing disc $D$. If $K''$ is isotopic to $K'$ then either $\partial D$ bounds a disc in the complement
of $K'$ or the crossing change occurs within one of $K'_1$, $K'_2$.
\end{theo}

\proof For an ordinary crossing the result is given as Theorem 2.1 in \cite{kn:to}.
The proof given in there works  for generalized crossings. \qed
\vskip 0.08in

The next lemma handles possibility (a) above as it reduces Theorem \ref{theo:last}
to the case that $K'$ is a prime knot.

\begin{lem} \label{lem:composite'} Let $K,K'$ be knots
such that $K\buildrel {(m, L_1)} \over \longrightarrow K'$,
where $L_1$ is a crossing circle for $K$. Suppose that
$K':=K'_1\# K'_2$ is a composite knot. Then, either
$L_1$ bounds a disc in the complement of $K$ or
$K$ is a connect sum $K = K_1\# K_2$ and there exist
$J\in \{K_1, K_2\}$ and $J'\in \{K'_1, K'_2\}$ such that
$J\buildrel {(m, L_1)} \over \longrightarrow J'$.
\end{lem}

\proof  By assumption there is an integer $r\neq 0$
so that the knot $K''$ obtained from $K'$ by a generalized crossing change of order $r$
is isotopic to $K'$. By Theorem \ref{theo:composite}, either $L_1$ bounds a disc in the complement
of $K'$ or the crossing change occurs on one of $K'_1$, $K'_2$; say on $K'_1$.
Thus, in particular, in the latter case $L_1$ is a crossing link for $K'_1$. Since
$K$ is obtained from $K'$ by twisting along $L_1$, $K$ is a, not necessarily non-trivial, connect sum
of the form $K_1\# K'_2$. By the uniqueness of knot decompositions it follows
that $K_1\buildrel {(m, L_1)} \over \longrightarrow K'_1$. \qed
\smallskip

\subsection {Dehn surgeries that create essential tori}
Let $M$ be a compact orientable 3-manifold.
For a collection
$\mathcal T$ of disjointly embedded, pairwise non-parallel,
essential tori
in $M$ we will use
$|\mathcal T|$ to denote 
the number of components of $\mathcal T$. By  Haken's finiteness theorem
(\cite{kn:hempel}, Lemma 13.2), the number
$$\tau(M)={\rm max}\, \{\,|{\mathcal T}|\, |\, {\mathcal T}\,
{\text{ a collection of tori as above}}\,\}$$
is well defined. A collection $\mathcal T$ for which $\tau(M)=|{\mathcal T}|$ will be called a {\em Haken system}.

In this section we will study
the behavior of essential tori under the various Dehn fillings from 
$N:=S^3\setminus \eta (K\cup L_1)$ to $N':=S^3\setminus \eta(K')$.
Since  
$N'$ is obtained from
$N$ by Dehn filling along 
$T_1:=\partial \eta(L_1)$, essential tori
in $N'$  occur in the following two ways:
\smallskip

Type $I$: An essential torus $T'\subset N'$ that can be
isotoped in $N\subset N'$; thus such a torus is the image
of an essential torus $T\subset N$.
\smallskip

Type $II$: An essential torus  $T'\subset N'$ that is the image
of an essential punctured torus $(P, \partial P)\subset (N, T_1)$,
such that each component of $\partial P$ is parallel on $T_1$ to the 
curve along which the Dehn filling from $N$ to $N'$ is done.
\smallskip

We begin with the following lemma that examines circumstances under which
twisting a knot that is geometrically essential inside a knotted solid torus $V$ yields
a knot that is geometrically inessential inside  $V$. In the notation of Definition \ref{defi:curve},
the lemma
implies that
an essential torus in $N$ either remains essential in
$N(s_{J})$, for all $\emptyset \neq J\subset \{1,\ldots, m\}$,
or it becomes inessential in all $N(s_{J})$.

\begin{lem} \label{lem:subcase1} Let $V\subset S^3$ be a knotted solid torus
and let $K_1\subset V$ be a knot that is geometrically essential in $V$. Let $D\subset {\rm int}(V)$ be a crossing disc
for $K_1$ and let $K_2$ be a knot obtained from $K_1$ by a non-trivial twist
along $D$. Suppose that $K_1$ is isotopic to $K_2$ in $S^3$. Then,
$K_2$ is geometrically essential in $V$. Furthermore, if $K_1$
is not the core of $V$ then $K_2$ is not the core of $V$.
\end{lem} 

\proof Suppose that $K_2$ is not geometrically essential in $V$.
Then there is
an embedded 3-ball $B\subset {\rm int}(V)$ that contains $K_2$. Since  making crossing changes
on $K_2$ doesn't change the homology class it represents in $V$, the winding number of $K_1$
in $V$ must be zero. Set $L:=\partial D$ and  $N:=S^3\setminus \eta(K\cup L)$.
Let $S$ be a Seifert surface for $K_1$ such that among all the surfaces bounded by $K_1$ in $N$,
$S$ has minimum genus. As usual we isotope $S$ so that $S\cap D$ is an arc $\alpha$
properly embedded on $S$. As in the proof of Theorem \ref{theo:genus},
$S$ gives rise to Seifert surfaces $S_1,S_2$ of $K_1, K_2$, respectively.
Now $K_1$ can be recovered from $K_2$ by twisting $\partial S_2$
along $\alpha$.

{\it Claim:} $L$ can be isotoped inside $B$ in the complement of $K_1$.
\vskip 0.04in

Since $K_1$ is obtained from $K_2$ by a
generalized crossing change supported on $L$ it follows that $K_1$  
lies in $B$.
Since this contradicts our assumption that $K_1$ is geometrically essential in $V$,
$K_2$ must be geometrically  essential in $V$. To finish the proof of the lemma, assuming the claim, observe that
that if $K_1$ is not the core $C$ of $V$, then $C$ is a companion knot of $K_1$. If
$K_2$ is the core of $V$, $C$ and $K_1$ are isotopic in $S^3$ which by
Schubert (\cite{kn:Knoten}) is impossible.

{\it Proof of Claim:} Since $K_1$, $K_2$ are isotopic in $S^3$ using Corollary
2.4 of \cite{kn:ga}, as used in the proof of Theorem \ref{theo:genus},
we see that $S_1$ (resp. $S_2$) is a minimum
genus surface for $K_1$ (resp. $K_2$) in $S^3$.
By assumption
$\partial V$ is a non-trivial companion torus of $K_1$. Since the winding number of $K_1$
in $V$ is zero, the intersection $S_1\cap \partial V$
(resp. $S_1\cap \partial V$) is homologically trivial
in $\partial V$. Thus we may replace the components
of $S_1\cap {\overline {S^3\setminus
V}}$, (resp. $S_2\cap {\overline {S^3\setminus
V}}$) with boundary parallel
annuli in ${\rm int}(V)$
to obtain a Seifert surface $S'_1$ (resp. $S'_2$)
inside  $V$. It follows, that
$S_1\cap {\overline {S^3\setminus
V}}$, (resp. $S_2\cap {\overline {S^3\setminus
V}}$) is a collection of annuli and 
$S'_1$ (resp. $S'_2$)
is a minimum genus Seifert surface for $K_1$ (resp. $K_2$).
Now $S'_2$ is a
minimum genus Seifert surface for $K_2$
such that $\alpha \subset S'_2$.
By assumption, $K_2$ lies inside $B$.
Since $S'_2$ is incompressible and $V$ is irreducible,
$S'_2$ can be isotoped in $B$ by a
sequence of disc trading isotopies in ${\rm int}(V)$.
But this isotopy will also bring $\alpha$ inside $B$
and thus $L$. \qed
\smallskip

Next we focus on the case that $N'$ is toroidal
and examine the circumstances under which $N'$ contains type  $II$ tori.
We have the following:

\begin{pro} \label{pro:torus} Let $K,K'$ be knots 
such that $K'$ is a non-trivial satellite but not composite. Suppose
that $K\buildrel {(m, L_1)} \over \longrightarrow K'$, where
$L_1$ is a crossing circle for $K$ and let
the notation be as in Definition \ref{defi:curve}. Then, at least
one of the following is true:
\vskip 0.04in

a) $L_1$ bounds an embedded disc in the complement
of $K$.
\vskip 0.04in

b) For every $\emptyset \neq J\subset \{1,\ldots, m\}$, $N(s_{J})$
has a Haken system that doesn't contain
tori of type $II$.
\vskip 0.04in

c) We have $m\leq 6$.
\end{pro}

\proof For $s\in {\bf Z}$, let $N(s)$
denote the 3-manifold obtained from $N$ by Dehn filling along $T_1$ with slope ${\textstyle {1\over s}}$.
Assume that $L_1$ doesn't bound an embedded disc
in the complement
of $K$ and that, for some $\emptyset \neq J_1\subset \{1,\ldots, m\}$,
$N(s_{J_1})$ admits a Haken system that contains tori of type $II$.
We claim that, for every $\emptyset \neq J\subset \{1,\ldots, m\}$,
$N(s_{J})$  has such a Haken system. To see this, first assume that
$N$ doesn't contain essential embedded tori. Then, since $N'=N(s_{J})$
and $K'$ is a non-trivial satellite the conclusion follows.
Suppose that $N$ contains essential embedded tori. By
Lemma \ref{lem:subcase1} it follows 
that an essential torus in $N$ either remains essential in
$N(s_{J})$, for all $\emptyset \neq J\subset \{1,\ldots, m\}$,
or it becomes inessential in all $N(s_{J})$ as above. Thus the number
of type $I$ tori in a Haken system of $N(s_{J})$ is the same for all
$J$ as above.
Thus, since we assume that $N(s_{J_1})$ has a Haken system containing tori of type $II$,
a Haken system of
$N(s_{J})$ must contain  tori of type $II$,
for every $\emptyset \neq J\subset \{1,\ldots, m\}$. We distinguish two cases:
\smallskip

{\it Case 1:} Suppose that $s_1, \ldots, s_m>0$ or  $s_1, \ldots, s_m<0$.
Let $s:= \sum_{j=1}^m s_j$ and recall that we assumed that $N$ is irreducible.
By our discussion above,
both of $N(s_1), N(s)$ contain essential embedded tori of type $II$.
By Theorem 1.1 of \cite{kn:go}, we must have
$$\Delta(s, s_1)\leq 5,\eqno (4.1)$$
\noindent where $\Delta(s, s_1)$ denotes the geometric intersection
on $T_1$
of the slopes represented by ${\textstyle {1\over {s_1}}}$, and ${\textstyle {1\over s}}$.
Since $\Delta(s, s_1)=|\sum_{j=2}^m s_j|$, and $|s_j|\geq 1$, in order for
(4.1) to be true we must have $m-1\leq 5$
or $m\leq 6$.
\smallskip

{\it Case 2:} Suppose that not all of
$s_1, \ldots, s_m$ have the same sign.
Suppose, without loss of generality, that
$s_1, \ldots, s_k>0$ and
$s_{k+1}, \ldots, s_m<0$.
Let $s:= \sum_{j=1}^k s_j$ and
$t:= \sum_{j=k+1}^m s_j$.
Since both of $N(s), N(t)$ contain essential embedded tori of type $II$,
by Theorem 1.1 of \cite{kn:go}
$$\Delta(s, t)\leq 5. \eqno(4.2)$$
But $\Delta(t, s)=s-t=\sum_{j=1}^m |s_j| $. Thus, in order for
(4.2) to be true we must have $m \leq 5$
and the result follows. \qed
\smallskip

\smallskip

\smallskip

Proposition \ref{pro:torus} and Lemma \ref{lem:composite'}  yield the following corollary:

\begin{corol} \label{corol:torus'} Let $K,K'$ be knots
and let $L_1$ be a crossing circle for $K$.
Suppose that the 3-manifold $N$
contains no essential embedded torus and that 
$K\buildrel {(m, L_1)} \over \longrightarrow K'$. If $K'$
is a non-trivial satellite, then
either $m\leq 6$ or $L_1$ bounds an embedded disc in the complement
of $K$.
\end{corol}
\smallskip

\subsection{ Hyperbolic and Seifert fibered manifolds} In this section we will deal with
the case that the manifold $N$ is atoroidal. As already mentioned, by Thurston's uniformization theorem
for Haken manifolds
(\cite{kn:th1}), $N$ is  either hyperbolic or a Seifert fibered manifold.

First we recall some  terminology about hyperbolic 3-manifolds.
Let $N$ be a hyperbolic 3-manifold with boundary and let $T_1$ a component of $\partial N$.
In ${\rm int}(N)$ there is a cusp, which is homeomorphic to $T_1\times [1, \infty)$,  associated with the torus
$T_1$. The cusp lifts to an infinite set, say $\mathcal H$,  of disjoint horoballs in
the hyperbolic space ${\bf H}^3$
which can be expanded so that each horoball in $\mathcal H$ has a point of tangency with some other.
The image of these horoballs under the projection ${\bf H}^3\longrightarrow {\rm int}(N)$,
is the {\em maximal horoball neighborhood} of $T_1$. The boundary ${\bf R^2}$ of each  horoball in ${\mathcal H}$
inherits a Euclidean metric from  ${\bf H}^3$ which in turn induces a Euclidean metric on $T_1$.
A slope $\bf s$ on $T_1$ defines a
primitive element in $\pi_1(T_1)$ corresponding to a Euclidean translation in 
${\bf R^2}$. The length of $\bf s$, denoted by $l({\mathbf s})$, is given by the length of corresponding
translation vector. 

Given a slope ${\mathbf s}$ on $T_1$, let us  use $N[{\mathbf s}]$ to denote the manifold obtained from
$N$ by Dehn filling along $T_1$ with slope ${\mathbf s}$. We remind the reader that in the case that the slope
${\mathbf s}$ is represented by ${\textstyle {1\over s}}$, for some $s\in {\bf Z}$, we use the notation $N(s)$
instead. Next we recall a result  of Cooper and Lackenby the proof of which relies on work of Thurston and Gromov.
We only state the result in the special case
needed here:

\begin{theo} \label{theo:cola}\ {\rm (Cooper-Lackenby,\ \cite{kn:cla})} \  Let $N'$ be a compact orientable manifold, with $\partial N'$
a collection of tori. Let $N$ be a hyperbolic manifold and let $\bf s$ be a slope on a toral component $T_1$
of $\partial N$ such that $N[{\mathbf s}]$ is
homeomorphic to $N'$. Suppose that the length of $\bf s$ on the maximal horoball of $T_1$ in ${\rm int}(N)$
is at least $2\pi+\epsilon$, for some $\epsilon > 0$. Then, for any given $N'$ and
$\epsilon > 0$, there is only a finite number of possibilities (up
to isometry) for $N$ and ${\mathbf s}$.
\end{theo}

\begin{remark} {\rm With the notation of Theorem \ref{theo:cola},
let $E$ denote the set of all slopes $\bf s$ on $T_1$,
such that $l({\mathbf s}) \leq 2\pi$.
It is a consequence of the Gromov-Thurston
``$2\pi$ " theorem that $E$ is finite. More specifically, the Gromov-Thurston
theorem (a proof of which is found in \cite{kn:bho}) states that
if $l({\mathbf s}) > 2\pi$, then $N[{\mathbf s}]$ admits a negatively curved metric. But in Theorem 11
of \cite{kn:bho}, Bleiler and Hodgson show that there can be at most 48
slopes on $T_1$ for which  $N[{\mathbf s}]$ admits no negatively curved metric.
Thus, there can be at most 48 slopes on $T_1$ with length $\leq 2\pi$.}
\end{remark}
\smallskip

Using Theorem \ref{theo:cola} we will prove the following
proposition which is a special case of Theorem \ref{theo:last} (compare possibility (b)
of \S 4.1):

\begin{pro} \label{pro:hyperbolic}Let $K,K'$ be knots
such 
that $K\buildrel {(m, L_1)} \over \longrightarrow K'$, where
$L_1$ is a crossing circle for $K$ and $m>0$.
Suppose that $N:=S^3\setminus \eta (K\cup L_1)$ is a hyperbolic manifold.
Then, there is a constant $b(K, \ K')$, that depends only on $K, K'$,
such that $m\leq b(K, \ K')$.
\end{pro}

\proof We will apply Theorem \ref{theo:cola} for the manifolds
$N:=S^3\setminus \eta (K\cup L_1)$, $N':=S^3\setminus \eta(K')$
and the component $T_1:=\partial \eta(L_1)$ of $\partial N$.
Let $s_1,\ldots, s_m$ be integers that satisfy Definition \ref{defi:curve}.
That is, for every $\emptyset \neq J\subset \{1,\ldots, m\}$,
$N(s_{J})$ is homeomorphic to $N'$.
By abusing the notation, for $r\in {\bf Z}$ we will use $l(r)$
to denote the length on $T_1$ of the slope represented by ${\textstyle {1\over r}}$.
Also, as in the proof of Proposition \ref{pro:torus}, we will use
$\Delta(r, t)$ to denote the geometric intersection on $T_1$
of the slopes represented by ${\textstyle {1\over r}}$, ${\textstyle {1\over t}}$.
Let $A(r, t)$ denote the area of the parallelogram in ${\bf R}^2$
spanned by the lifts of these slopes and let $A(T_1)$ denote the area
of a fundamental domain of the torus $T_1$. It is known that
$A(T_1) \geq {\textstyle {\sqrt 3\over 2}}$ (see, \cite{kn:bho})
and that $\Delta(r, t)$
is the quotient of $A(r, t)$ by $A(T_1)$. Thus, for every $r,t\in {\bf Z}$, we have

$$l(r) l(t) \geq \Delta(r, t){\sqrt 3\over 2}. \eqno (4.3)$$

Let $\lambda>0 $ denote the length of a meridian of $T_1$; in fact it is known that $\lambda\geq 1$.
Assume on the contrary that no constant $b(K, \ K')$ as in the statement
of the proposition exists. Then, there exist infinitely many
integers $s$ such that $N(s)$ is homeomorphic to $N'$. Applying (4.3)
for $l(s)$ and $\lambda$ we obtain
$$l(s)\geq |s| {{\sqrt 3\over 2\lambda}}.$$ 
Thus, for $|s|\geq {\textstyle {{{4\pi \lambda +2\lambda}}\over \sqrt 3}}$
we have  $l(s)\geq 2\pi +1$.
But then for $\epsilon=1$,
we have infinitely many integers such that $l(s)\geq 2\pi +\epsilon$
and $N(s)$ is homeomorphic to $N'$. Since this contradicts
Theorem \ref{theo:cola} the proof of the Proposition is finished. \qed

\smallskip

\smallskip

Next we turn our attention to the case where
$N:=S^3\setminus \eta (K\cup L_1)$ is an atoroidal Seifert fibered space. Since $N$ is embedded in $S^3$ it is orientable.
It is know that an orientable, atoroidal
Seifert fibered space with two boundary components is
either
a {\em cable space} or a trivial torus bundle $T^2\times I$.
Let us recall how a cable space is formed: Let $V''\subset V'\subset S^3$ be concentric
solid tori. Let $J$ be a simple closed curve on $\partial V''$ having slope
${\textstyle {a\over b}}$, for some $a,b\in {\bf Z}$ with $|b|\geq 2$.
The complement $X:=V'\setminus {\rm int}(\eta(J))$ is a ${\textstyle {a\over b}}$-cable space.
Topologically, $X$ is a Seifert fibered space over the annulus with one exceptional
fiber of multiplicity $|b|$.
We show the following:

\begin{lem} \label{lem:Seifert}Let $K,K'$ be knots 
such 
that $K\buildrel {(m, L_1)} \over \longrightarrow K'$, where
$L_1$ is a crossing circle for $K$ and $m>0$.
Suppose that $N:=S^3\setminus \eta (K\cup L_1)$ is an irreducible,  atoroidal Seifert fibered space.
Then, there is a constant $b(K, \ K')$
such that $m\leq b(K, \ K')$.
\end{lem} 

\proof As discussed above, $N$ is either a {cable space} or a torus bundle $T^2\times I$.
Note, however, that in a cable space the cores of the solid tori bounded in $S^3$ by the two components
of $\partial N$  have non-zero linking number. Thus, since the linking number of $K$
and $L_1$ is zero, $N$ cannot be a cable space.
Hence, 
we only have to consider the case where
$N \cong T^2\times I$. Suppose $T_1=T^2\times \{ 1\}$ and
$T_2:=\partial \eta(K)= T^2\times \{ 0 \}$. By assumption there is a slope $\bf s$
on $T_1$ such that the Dehn filling of $T_1$ along ${\mathbf s}$ produces $N'$. Now
${\bf s}$ corresponds to a simple closed curve
on $T_2$ that must compress in $N'$. By Dehn's Lemma, $K'$ must be the unknot.
It follows that either $g(K)>g(K')$ or $K$ is the unknot. In the later case,
we obtain that $L_1$ bounds a disc disjoint from $K$
contrary to our assumption that $N$ is irreducible. Thus, $g(K)>g(K')$
and the conclusion follows from Theorem \ref{theo:lower}. \qed.
\smallskip

The following Proposition
complements nicely Corollary \ref{corol:torus'}. We point out that the proposition is not needed
for the proof of the main result.
Hence
a reader 
eager to get to the proof
of Theorem \ref{theo:last} can
move to the next section without loss of continuity.

\begin{pro} \label{pro:hyperbolic'} Let $K,K'$ 
be non-isotopic hyperbolic knots. Suppose 
there exists  a
crossing circle $L_1$ for $K$  such that
$K\buildrel {(m, L_1)} \over \longrightarrow K'$, for some $m\geq 6$.
Then, for given $K$ and $K'$, there is only a finite number of possibilities for $m$ and for
$L_1$ up to isotopy in the complement of $K$.

\end{pro}
\proof As before, let $N:=S^3\setminus \eta (K\cup L_1)$, $N':=S^3\setminus \eta(K')$
and let $D$ be a crossing disc for
$L_1$. Since $K$ is not isotopic to $K'$, $N$ is irreducible
and $\partial$-irreducible.

 {\it Claim.} $N$ is atoroidal.

 {\it Proof of Claim.} Suppose that $N$ contains an embedded essential torus $T$
and let $V$ denote the solid torus bounded by $T$ in $S^3$.
If $L_1$ cannot be isotoped to lie in ${\rm int}V$ then $D\cap T$ contains
a component whose interior in $D$ is pierced exactly once by $K$.
This implies that $T$ is parallel to $\partial \eta(K)$ in $N$; a contradiction.
Thus, $L_1$ can be isotoped to lie inside $V$. Now let $S$ be a Seifert surface of $K$
that is taut in $N$. After isotopy, $D\cap S$
is an arc $\alpha$ that is essential on $S$. By Theorem \ref{theo:genus},
$S$ remains of minimum genus in at least one of $N'':=S^3\setminus \eta(K)$, $N'$.
Assume $S$ remains of minimum genus in  $N'$; the other case is completely analogous.
Since $K,K'$ are hyperbolic
$T$ becomes inessential in both of $N''$, $N'$.
But since $K,K'$ are related by a generalized crossing change, either
$T$ becomes boundary parallel in both of $N''$, $N'$ or
it becomes compressible in both of them. 
First suppose that $T$ is  boundary parallel in both of $N''$, $N'$:
Then it follows that the arc $\alpha$ is inessential on $S$ and $K$ is isotopic to $K'$; a contradiction.
Now suppose that $T$ is compressible in both of $N'',N'$: Then,
both of $K, K'$ are inessential in $V$ and they can be isotoped to
lie in a 3-ball $B\subset {\rm int}V$. By an argument similar to this in the proof
of Lemma \ref{lem:subcase1} we can conclude that $\alpha$, and thus $L_1$, can be isotoped to lie in $B$.
But this contradicts the assumption that $T$ is essential in $N$ and finishes the proof of the claim.
\vskip 0.04in

\smallskip

To continue observe that the argument of the proof of Lemma \ref{lem:Seifert} shows that
if $N$ is  a Seifert fibered space then $K'$
is the unknot. But this is impossible since we assumed that $K'$ is hyperbolic.
Thus, by \cite{kn:th1},
$N$ is hyperbolic. Let $s_1,\ldots, s_m$ be integers that satisfy Definition \ref{defi:curve}
for $K,K'$. Thus we have $2^m-1$ integers $s$, with $N(s)=N'$.
Now \cite{kn:bho} implies that we can have at most 48 integers so that the corresponding slopes have lengths $\leq 2\pi$ on $T_1$.
Since $m\geq 6$ we have $2^m-1>48$. Thus we have $k_m:=2^m-49>0$
integers $s$ such that $l(s)>2\pi$ and $N(s)=N'$. By Theorem \ref{theo:cola},
there is only a finite number of possibilities 
(up
to isometry) for $N$ and ${s}$.
Now the conclusion follows. \qed

\begin{remark} {\rm  Proposition \ref{pro:hyperbolic'}
implies Theorem \ref{theo:last}, and thus Theorem \ref{theo:main},
in the case that $K$, $K'$ are hyperbolic.}
\end{remark}

\subsection{The proof of Theorem \ref{theo:last}} In this subsection we give the proof
of Theorem \ref{theo:last}. We will need the following theorem which is a special case of a result of 
McCullough proven in \cite{kn:appen}.

\begin{theo} {\rm (McCullough,\ \cite{kn:appen})} Let $M$ be a compact orientable $3$-manifold, and let
$C$ be a simple loop in $\partial
M$. Suppose that $h\colon M\to M$ is a homeomorphism whose restriction to
$\partial M$ is isotopic to a nontrivial power of a Dehn twist about
$C$. Then, $C$ bounds a disc in $M$.
\label{thm:mainthm}
\end{theo}

\smallskip

Before we embark on the proof of Theorem \ref{theo:last}, 
we recall that for a compact orientable 3-manifold $M$,
$\tau(M)$ denotes the cardinality of a Haken system of tori (see subsection \S 4.3).
In particular, $M$ is atoroidal if and only if $\tau(M)=0$.
\smallskip
\smallskip

\proof {[Proof of Theorem \ref{theo:last}]}
Let $K, K'$ be knots and let $L_1$ be a
crossing circle of $K$ such that
$K\buildrel {(m, L_1)} \over \longrightarrow K'$.
As before we set
$N:=S^3\setminus \eta (K\cup L_1)$
and $N':=S^3\setminus \eta (K')$.
If
$g(K)> g(K')$, by Theorem \ref{theo:lower},
we have $m\leq 3g(K)-1$. Thus, in this case,
we can take
$b(K,\,K'):=3g(K)-1$ and Theorem \ref{theo:last} holds.
Hence, we only have to consider that case
that $g(K)\leq g(K')$.

Next  we consider the complexity
$$\rho=\rho(K,\,K',\,L_1): =\tau(N).$$

First, suppose that $\rho=0$, that is $N$ is atoroidal. Then,
$N$ is either hyperbolic or a Seifert fibered manifold (\cite{kn:th1}). In the former case,
the conclusion of the theorem follows from Proposition \ref{pro:hyperbolic};
in the later case it follows from Lemma \ref{lem:Seifert}.

Assume now that $\tau(N)>0$; that is $N$ is toroidal.
Suppose, inductively, that for every triple
$K_1, K_1', L_1'$, with $\rho(K_1,\, K'_1,\,L'_1)<r$, 
there is a constant $d=d(K_1,\,K_1')$
such that: If $K_1\buildrel {(m, L_1')} \over \longrightarrow K_1'$,
then either 
$m\leq d$ or
$L_1'$ bounds an embedded disc in the complement
of $K_1$.
Let $K, K', L_1$ be knots and a crossing circle for $K$, such that
$K\buildrel {(m, L_1)} \over \longrightarrow K'$
and $\rho(K,\,K',\,L_1)=r$.
Let $s_1,\ldots, s_m$ be integers satisfying Definition \ref{defi:curve}
for $K, K'$ and $L_1$.
For every $\emptyset \neq J\subset \{1,\ldots, m\}$,
let $N(s_{J})$ be the 3-manifold obtained from
$N$ by Dehn filling of $\partial \eta(L_1)$ with slope
${\textstyle {1\over s_J}}$.
By assumption, $N'=N(s_{J})$.
Assume, for a moment, that
for some $\emptyset \neq J_1\subset \{1,\ldots, m\}$,
$N(s_{J_1})$ contains essential embedded tori of type $II$.
Then 
Proposition \ref{pro:torus} implies that either
$m\leq 6$ or
$L_1$ bounds an embedded disc in the complement
of $K$. Hence, in this case, the conclusion
of the theorem is true for $K, K', L_1$, with
$b(K,\ K'):=6$. Thus we may assume that,
for every $\emptyset \neq J\subset \{1,\ldots, m\}$,
$N(s_{J})$ 
doesn't contain essential embedded tori of type $II$.
\smallskip

We will show the following:

{\it Claim 1:} There exist knots $K_1, K_1'$ and a crossing circle $L_1'$ for
$K_1$ such that:

(1)  $K_1\buildrel {(m, L_1')} \over \longrightarrow K_1'$
and $\rho(K_1,\,K'_1,\,L'_1)< \rho(K,\, K',\,L_1)=r$. 

(2) If $L'_1$ bounds an embedded disc in the complement of $K_1$
then $L_1$ bounds an embedded disc in the complement of $K$.

{\it The proof of the theorem assuming Claim 1:}
By induction, 
there is $d=d(K_1, \ K_1')$
such that either $m\leq d$ or 
$L'_1$ bounds a disc in the complement of $K_1$.
Let ${\mathcal K}_m$ denote the set of all pairs of knots $K_1, K'_1$
such that there exists a crossing circle $L_1'$ for
$K_1$ satisfying properties (1) and (2) of Claim 1.
Define

$$b=b(K,\,K'):={\rm min}\,\{\,d(K_1,\,K_1')\,|\,K_1, K'_1\in{\mathcal K}_m\, \}.$$
Clearly $b$ satisfies the conclusion of the statement of the theorem.
\smallskip

{\it Proof of Claim 1:} Let $T$ be an  essential embedded torus
in $N$. Since $T$ is essential in $N$, $T$ has to be knotted.
Let $V$ denote the solid torus component of 
$S^3\setminus T$. Note that $K$ must
lie inside $V$. For, otherwise
$L_1$ must be geometrically essential in $V$ and thus it can't be the unknot.
There are various cases to consider
according to whether $L_1$ lies outside or inside $V$.
\vskip 0.06in

{\it Case 1:} Suppose that $L_1$ lies outside $V$
and it cannot be isotoped to lie inside $V$.
Now $K$ is a non-trivial satellite with companion torus $T$.
Let $D_1$ be a crossing disc bounded by $L_1$.
Notice that if all the components of $D_1\cap T$
were either homotopically trivial in $D_1\setminus(D_1\cap K)$ 
or parallel to $\partial D_1$,
then we would be able to isotope $L_1$
inside $V$ contrary to our assumption.
Thus $D_1\cap T$
contains a component that encircles a
single point of the intersection $K\cap D_1$.
This implies that the winding number of $K$ in $V$ is
one. Since $T$ is  essential 
in $N$ we conclude that $K$ is composite, say $K:={K}_1\# {K}_2$,
and $T$ is the follow-swallow torus. Moreover,
the generalized crossings realized by the surgeries on $L_1$ occur
along a summand of $K$, say along ${K}_1$. By the uniqueness of prime decompositions
of knots, it follows
that there exists a (not necessarily non-trivial) knot $K'_1$,
such that
$K'={ K'}_1\# {K}_2$ and
$K_1\buildrel {(m, L_1)} \over \longrightarrow K'_1$.
Set $N_1:=S^3\setminus \eta (K_1\cup L_1)$
and $N'_1:=S^3\setminus \eta (K'_1)$.
Clearly, 
$\tau(N_1)<\tau(N)$.
Thus,
$\rho(K_1,\,K'_1,\,L_1)< \rho(K,\,K',\,L_1)$
and part (1) of the claim
has been proven in this case. To see part (2) notice that
if $L_1$ bounds a disc $D$ in the complement
of $K_1$, we may assume $D\cap K=\emptyset$.
\vskip 0.06in

{\it Case 2:} Suppose that $L_1$ can be isotoped to lie inside $V$.
Now the link $K\cup L_1$ is a non-trivial satellite with companion torus
$T$. We can find a standardly
embedded solid torus $V_1\subset S^3$, and
a 2-component link $(K_1\cup L'_1) \subset V_1$ such that:
i) $K_1\cup L_1'$ is geometrically essential in $V_1$; 
ii) $L'_1$ is a crossing disc for $K_1$; and iii)
there is 
a homeomorphism  $f: V_1 \longrightarrow V$
such that $f(K_1)=K$ and $f(L'_1)=L_1$ and $f$ preserves the longitudes
of $V_1$ and $V$. In other words,
$K_1\cup L'_1$ is the model link for the satellite.
Let $\mathcal T$ be a Haken system for $N$ containing $T$.
We will assume that the torus $T$
is innermost; i.e. the
boundary of the component of $N\setminus {\mathcal T}$ that contains $T$ also contains
$\partial \eta(K)$. By twisting along $L_1$ if necessary, we may without loss of generality assume that
${\bar V}:=\overline{V\setminus ( K\cup L_1)}$
is atoroidal. Then, ${\bar V}_1:=\overline{V_1\setminus ( K_1\cup L'_1)}$
is also atoroidal.
For every $\emptyset \neq J\subset \{1,\ldots, m\}$,
let $K(s_J)$ denote the knot obtained from $K_1$ by performing
${\textstyle {1\over s_J}}$-surgery on
$L'_1$. By assumption the knots $f(K(s_J))$ are all isotopic
to $K'$.
\vskip 0.06in

{\it Subcase 1:} There is $\emptyset \neq J_1\subset \{1,\ldots, m\}$, such that
$\partial V$ is compressible in $V\setminus f(K(s_{J_1}))$.
By Lemma \ref{lem:subcase1},
for every $\emptyset \neq J
\subset \{1,\ldots, m\}$,
$\partial V$ is compressible in $V\setminus f(K(s_J))$. It follows that there is an
embedded 3-ball $B\subset {\rm int}(V)$ such that: i) $f(K(s_J))\subset {\rm int }(B)$, for
 every $\emptyset \neq J
\subset \{1,\ldots, m\}$; and ii) the isotopy from $f(K(s_{J_1}))$
to $f(K(s_{J_2}))$ can be realized inside $B$, for every $J_1\neq J_2$ as above.
From this observation it follows that there
is a knot $K_1'\subset {\rm int}(V_1)$
such that $f(K_1')=K'$ and
$K_1\buildrel {(m, L_1')} \over \longrightarrow K'_1$
in $V_1$. Let $N_1:=S^3\setminus \eta(K_1\cup L'_1)$
and $N'_1:=S^3\setminus \eta (K'_1)$.
Clearly,
$\tau(N_1)<\tau(N)$.
Hence, $\rho(K_1,\,K'_1,\,L_1)< \rho(K,\, K',\,L_1)$
and the part (1) of Claim 1 has been proven.

We will prove part (2) of Claim 1 for this subcase together
with the next subcase.
\vskip 0.06in

{\it Subcase 2:} For every $\emptyset \neq J\subset \{1,\ldots, m\}$,
$f(K(s_{J}))$ is geometrically essential in $V$.
By Lemma \ref{lem:composite'}, the conclusion of the claim is true if
$K'$ is composite. Thus, we may assume that
$K'$ is a prime knot. 
In this case, we claim that, for every $\emptyset \neq J\subset \{1,\ldots, m\}$,
there is an orientation preserving
homeomorphism $\phi:S^3\longrightarrow S^3$ such that $\phi(V)=V$
and $\phi(f(K(s_{J_1})))=f(K(s_{J_2}))$. 
Since we assumed that
$N(s_{J_1})$, $N(s_{J_1})$ do not contain essential tori of type
$II$,
$T$ remains innermost in the complement of $f(K(s_{J_1}))$, $f(K(s_{J_2}))$.
By the uniqueness of the torus decomposition of knot 
complements \cite{kn:JS} or the uniqueness of satellite structures of knots
\cite{kn:Knoten}, there is an orientation preserving
homeomorphism $\phi:S^3\longrightarrow S^3$ such that $\phi(V)\cap V=\emptyset $
and ${\bar K}:=\phi(f(K(s_{J_1})))=f(K(s_{J_2}))$ (compare, Lemma 2.3 of  \cite{kn:Motegi}).
Since $T$ is innermost
in ${\bar V}$, we have
$S^3\setminus {\rm int}(V)\subset {\rm int}(\phi(S^3\setminus {\rm int}(V)))$ or
$\phi(S^3\setminus {\rm int}(V)) \subset {\rm int}(S^3\setminus {\rm int}(V))$.
In both cases, by Haken's finiteness theorem, it follows that
$T$ and $\phi(T)$ are parallel in the complement of  $\bar K$. Thus after an ambient isotopy, leaving
$\bar K$ fixed, 
we have $\phi(V)=V$.
Let $h=f\circ\phi\circ f^{-1}:V_1
\longrightarrow V_1$. Then $h$ preserves the longitude of $V_1$
up to a sign
and $h(K(s_{J_1}))=K(s_{J_2})$.
So, in particular, the knots $K(s_{J_1})$ and $K(s_{J_2})$
are isotopic in $S^3$.
Let $K_1'$ denote the knot type in $S^3$ of
$\{K(s_J)\}_{J\subset {\{1,\ldots, m\}}}$. By our earlier assumptions, 
$K_1\buildrel {(m, L_1')} \over \longrightarrow K'_1$.
Let $N_1:=S^3\setminus \eta(K_1\cup L'_1)$
and $N'_1:=S^3\setminus \eta (K'_1)$.
Clearly,
$\tau(N_1)<\tau(N)$.
Thus part (1) of Claim 1 has been proven also in this subcase.
\vskip 0.06in

We now prove part (2) of Claim 1 for both subcases.
Note that it is enough to show that if $L_1'$
bounds an embedded disc, say $D'$, 
in the complement of $K_1$ in $S^3$, then it bounds one
inside $V_1$. 

Let
$D_1'\subset V_1$ be a crossing disc bounded by
$L_1'$ and such that ${\rm int}(D')\cap {\rm int}(D_1')=\emptyset$.
Since $\partial V_1$ is incompressible in $V_1\setminus K_1$, after
a cut and paste argument, we may assume that 
$E=D_1'\cup(D\cap V_1)$ is a proper annulus whose boundary are
longitudes of $V_1$.

{\vspace{.03in}}
\begin{figure}[ht]
\centerline{\psfig{figure=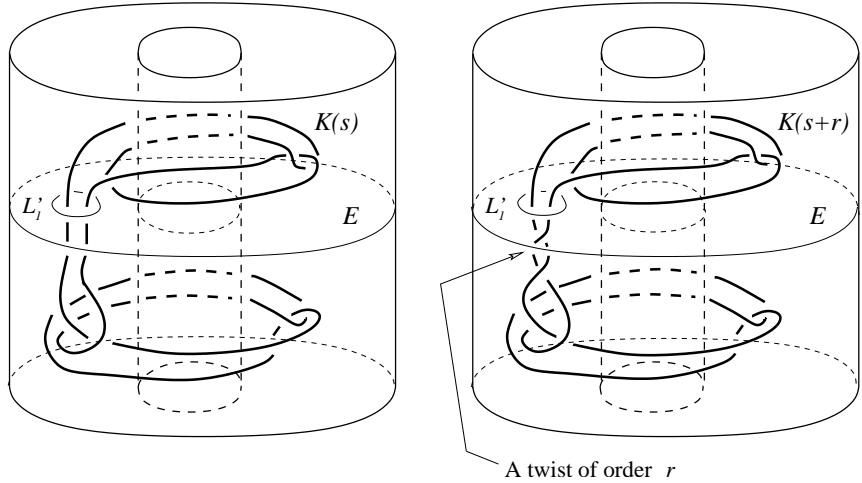,height=2.50in, clip=}}
\caption{The annulus $E$ contains the crossing circle $L'_1$ and
separates $V_1$ into solid tori $V'_1$ (part above $E$) and $V''_1$
(part below $E$.) 
In $V''_1$ the knots $K(s),\,K(s+r)$ differ by a twist of order $r$ along $D'_1$.}
\end{figure}

By assumption, in both subcases, there exist non-zero integers $s, r$,
such that
$K(s)$ and $K(s+r)$ are isotopic in $S^3$. 
Here, $K(s)$ and $K(s+r)$ 
denotes the knots obtained from $K_1$
by a twist along $L_1'$ of order $s$ and $s+r$ respectively.
Let $\hat h:S^3\longrightarrow S^3$ denote the extension 
of $h:V_1\longrightarrow V_1$ to $S^3$. We assume that $\hat h$ fixes the
core circle $C_1$ of the complementary solid torus of $V_1$. 
Since the 2-sphere $D\cup D_1'$ gives the same (possible trivial)
connected sum decomposition of $K_1'=K(s)=K(s+r)$ in $S^3$, we may 
assume that $\hat h(D)=D$ and $\hat h(D_1')=D_1'$ up to an isotopy.
During this isotopy of $\hat h$, $\hat h(C_1)$ and $\hat h(V_1)$
remain disjoint. So we may assume that at the end of the isotopy, we still
have $\hat h(V_1)=V_1$. Thus, we can assume that $h(E)=E$.

The annulus $E$ cuts $V_1$ into 
two solid tori $V_1'$ and $V_1''$. See Figure 4,
where the solid torus above $E$ is $V_1'$ and below $E$ is $V_1''$.  
We have either $h(V_1')=V_1'$ and $h(V_1'')=V_1''$ or
$h(V_1')=V_1''$ and $h(V_1'')=V_1'$. 
In the case when  $h(V_1')=V_1'$ and $h(V_1'')=V_1''$, we may assume that
$h|\partial V_1=\rm id$ and $h|E=\rm id$. Thus $K(s+r)\cap V_1'=K(s)\cap 
V_1'$ and $K(s+r)\cap V_1''$
is equal to $K(s)\cap V_1''$ twisted by a twist of order $r$ along $L_1'$.
Let $M$ denote the 3-manifold obtained from $V_1''\setminus(V_1''\cap K(s))$
by attaching a 2-handle to $\partial V_1''\cap E$ along  $K(s)\cap V_1''$.
Now $h|\partial M$ can be realized by
a Dehn twist of order $r$  along $L_1'$.
By Theorem \ref{thm:mainthm},
$L_1'$ must bound a disc in $M$. In order words, $L_1'$ bounds a disc in $V_1\setminus K(s)$.
This implies that $L_1'$ bounds a disc in $V_1\setminus K_1$ 

In the case when $h(V_1')=V_1''$ and $h(V_1'')=V_1'$, we may assume that
$h|\partial V_1$ and $h|E$ are  rotations of $180^\circ$ with an
axis on $E$ passing through the intersection points of $D_1'$ with
$K(s)$ and $K(s+r)$. Thus $K(s+r)\cap V_1'$ and $K(s)\cap 
V_1''$ differ by a rotation, and $K(s+r)\cap V_1''$ 
is equal to $K(s)\cap V_1'$ twisted by a twist of order $r$ along $L_1'$ 
followed by a rotation. Now we consider the 3-manifold
$N$ obtained from $V_1'\setminus(V_1'\cap K(s))$
by attaching a 2-handle to $\partial V_1'\cap E$ along  $K(s)\cap V_1'$.
As above we conclude that
a Dehn twist of order $r$ along 
$L_1'$  extends to 
$N$ and we complete the argument by
applying Theorem \ref{thm:mainthm}.  \qed 
\smallskip

\section{Applications and examples}
\subsection{Applications to nugatory crossings }
Recall that a crossing of a
knot $K$ with crossing disc $D$ is called {\it nugatory}
if $\partial D$ bounds a disc
disjoint from $K$. This disc and $D$ bound a
2-sphere that decomposes $K$ into a connected sum,
where some of the summands may be trivial.
Clearly, changing a nugatory crossing doesn't
change the isotopy class of a knot.
An outstanding open question is whether
the converse is true (see Problem 1.58 of Kirby's Problem List (\cite{kn:kirby}):

\begin{question} \label{question:kirby}\ {\rm ( Problem 1.58, \cite{kn:kirby})} \ If a crossing change in a knot $K$ yields a knot
isotopic to $K$ is the crossing
nugatory?
\end{question}

The answer is know to be {\em yes} in the case when $K$ is the unknot
(\cite{kn:st}) and when $K$ is a
2-bridge knot (\cite{kn:to}). In \cite{kn:to}, I. Torisu conjectures that the
answer is always {\em yes}. Our results in Section five
yield the following corollary that shows that an {\em essential}
crossing circle of a knot $K$  can admit at most finitely many twists that do not change
the isotopy type of $K$:

\begin{corol} \label{corol:strong}
For a crossing of a knot K, with crossing disc $D$, let $K(r)$ denote the knot obtained
by a twist of order $r$  along $D$. The crossing is nugatory
if and only if
$K(r)$ is isotopic to $K$ for all $r\in \Z$. 
\end{corol}
\proof One direction of the corollary is clear. To obtain the other direction apply Theorem
\ref{theo:last} for $K=K'$. \qed
\smallskip

In the view of Corollary \ref{corol:strong},  Question \ref{question:kirby}
is reduced to the following:
With the same setting as in Corollary \ref{corol:strong}, let $K_+:=K$ and
$K_-:=K(1)$. If $K_-$ is isotopic to $K_+$ is it true that $K(r)$ is isotopic to $K$, for all $r \in \Z$? 
\smallskip

\subsection{Examples}
In this subsection, we outline  some methods
that for every $n>0$ construct knots
$K, K'$ with \ad. 
It is known that given $n\in {\bf N}$
there exists a plethora of knots that are $n$-adjacent to the unknot.
In fact,
\cite{kn:ak} provides a method for constructing all such knots.
It is easy to see that given knots $K,K'$
such that $K_1$ is $n$-adjacent to the unknot,
the connected sum $K:=K_1\# K'$ is $n$-adjacent to $K'$. Clearly,
if $K_1$ is non-trivial then $g(K)>g(K')$.
To construct examples $K,K'$ in which $K$ is not
composite, at least in an obvious way, one can proceed as follows:
For $n>0$ let $K_1$ be a knot that is $n$-adjacent to the unknot
and let $V_1\subset S^3$ be a standard solid torus.
We can embed $K_1$ in $V_1$ so that i) it has non-zero winding number;
and ii) it is $n$-adjacent to the core of $V_1$ inside $V_1$.
Note that there might be many different ways of doing so.
Now let $f:V_1\longrightarrow S^3$ be any embedding that knots
$V_1$. Set $V:=f(V_1)$, $K:=f(K_1)$  and let $K'$ denote the core
of $V$. By construction, \ad.
Since $K_1$ has non-zero winding number in $V_1$ we have
$g(K)>g(K')$ (see, for example, \cite{kn:bz}).

We will say that two ordered pairs of knots
$(K_1, \ K_1')$, $(K_2, \ K_2')$ are isotopic iff
$K_1$ is isotopic to $K_2$ and   $K_1'$
is isotopic to $K_2'$.
From our discussion above we obtain the following:

\begin{pro} \label{pro:examples}
For every $n\in {\bf N}$ there exist infinitely many non-isotopic pairs
of knots $(K, \ K')$ such that \ad\ and $g(K)>g(K')$.
\end{pro}

\begin{remark} {\rm We should point out that, at this point, we don't know of any examples
of knots $(K, \ K')$ such that \ad\ and $g(K)<g(K')$. 
In fact the results 
of \cite{kn:k1}, and further examples constructed by Torisu \cite{kn:to1}, prompt the following question: Is it true that if \ad\, for some $n>1$, then either 
$g(K)>g(K')$ or $K$ is isotopic to $K'$?}
\end{remark}

\end{document}